\documentclass[12pt]{article} 
\usepackage{latexsym, 
amscd, 
amsthm, 
graphicx, epsfig, amssymb}

\newtheorem{thm}{Theorem}[section]
\newtheorem{defn}[thm]{Definition}

\newtheorem{lemma}[thm]{Lemma}

\newcommand{\halmos}{\rule{1ex}{1.4ex}}

\newcommand{\nn}{\nonumber \\}

 \newcommand{\pf}{{\it Proof.}\hspace{2ex}}
 
 \newcommand{\epfv}{\hspace*{\fill}\mbox{$\halmos$}\vspace{1em}}

\newcommand{\C}{\mathbb{C}}
\newcommand{\Z}{\mathbb{Z}}

\newcommand{\N}{\mathbb{N}}

\newcommand{\one}{\mathbf{1}}

\title{ {\bf First and second cohomologies 
of grading-restricted vertex algebras} }
\date{}
\author{Yi-Zhi Huang}
\begin{document}

\bibliographystyle{alpha}
\maketitle

\begin{abstract}
Let $V$ be a grading-restricted vertex algebra and $W$ a $V$-module.
We show that for any $m\in \Z_{+}$, the first cohomology $H^{1}_{m}(V, W)$ 
of $V$  with coefficients in $W$
introduced by the author is 
linearly isomorphic 
to the space of derivations from $V$ to $W$. In particular, 
$H^{1}_{m}(V, W)$ for 
$m\in \N$ are equal (and can be denoted using the same notation
$H^{1}(V, W)$). 
We also show that the 
second cohomology $H^{2}_{\frac{1}{2}}(V, W)$ of $V$  with coefficients in $W$ 
introduced by the author
corresponds bijectively to the set of equivalence classes of square-zero extensions of 
$V$ by $W$. In the case that $W=V$, we show that the 
second cohomology $H^{2}_{\frac{1}{2}}(V, V)$
corresponds bijectively to the set of equivalence classes of
first order deformations 
of $V$. 
\end{abstract}

\renewcommand{\theequation}{\thesection.\arabic{equation}}
\renewcommand{\thethm}{\thesection.\arabic{thm}}
\setcounter{equation}{0}
\setcounter{thm}{0}

\section{Introduction}

The present paper is a sequel to the paper \cite{H}. We discuss the 
first and second cohomologies of grading-restricted 
vertex algebras introduced by the author in that paper.

Let $V$ be a grading-restricted vertex algebra and $W$ a $V$-module.
Recall from \cite{H} that for each $m\in \Z_{+}$ and $n\in \N$, we have an
$n$-th cohomology $H^{n}_{m}(V, W)$ 
of $V$  with coefficients in $W$. For each $n\in \N$, We also have an
$n$-th cohomology $H^{n}_{\infty}(V, W)$ of $V$  
with coefficients in $W$
which is isomorphic to the 
inverse limit of the inverse system $\{H^{n}_{m}(V, W)\}_{m\in \Z_{+}}$. 
We also have
an additional second cohomology $H^{2}_{\frac{1}{2}}(V, W)$
of $V$  with coefficients in $W$.
In the present paper, we discuss only $H^{1}_{m}(V, W)$ for $m\in \Z_{+}$
and $H^{2}_{\frac{1}{2}}(V, W)$.

Let $V$ be a grading-restricted vertex algebra and $W$ a 
$V$-module. 
A grading-preserving linear map $f: V \to W$ is called a 
{\it derivation} if 
\begin{eqnarray*}
f(Y_{V}(u, z)v)&=&Y_{WV}^{W}(f(u), z)v+Y_{W}(u, z)f(v)\nn
&=&e^{zL(-1)}Y_{W}(v, -z)f(u)+Y_{W}(u, z)f(v)
\end{eqnarray*}
for $u, v\in V$.
We use $\mbox{\rm Der}\;(V, W)$ to denote the space of all such derivations.
We have the following result for the first cohomologies of $V$
with coefficients in $W$:

\begin{thm}\label{1st-coh}
Let $V$ be a grading-restricted vertex algebra and 
$W$ a $V$-module. Then  $H_{m}^{1}(V, W)$ 
is linearly isomorphic to the space of derivations from $V$ to $W$
for any $m\in \Z_{+}$, that is,
$H_{m}^{1}(V, W)$ 
is linearly isomorphic to $\mbox{\rm Der}\;(V, W)$ for any $m\in \Z_{+}$.
\end{thm}

In particular, 
$H^{1}_{m}(V, W)$ for 
$m\in \N$ are isomorphic (and can be denoted using the same notation
$H^{1}(V, W)$). 

\begin{defn}
{\rm Let $V$ be a grading-restricted vertex algebra.
A {\it square-zero ideal of $V$} is an ideal $W$ of $V$ such that 
for any $u, v\in W$, $Y_{V}(u, x)v=0$.}
\end{defn}

\begin{defn}
{\rm Let $V$ be a grading-restricted vertex algebra and $W$ 
a $\Z$-graded $V$-module. 
A {\it square-zero extension $(\Lambda, f, g)$ of $V$ by $W$}
is a grading-restricted vertex algebra $\Lambda$ 
together with a surjective homomorphism $f: \Lambda \to V$ 
of grading-restricted vertex algebras such that $\ker f$
is a square-zero ideal of $\Lambda$ (and therefore a $V$-module)
and an injective homomorphism $g$ of $V$-modules from $W$ to $\Lambda$ 
such that $g(W)=\ker f$. 
Two square-zero extensions $(\Lambda_{1}, f_{1}, g_{1})$ and 
$(\Lambda_{2}, f_{2}, g_{2})$ of $V$ by $W$ are {\it equivalent} if there 
exists an isomorphism of grading-restricted vertex algebras
$h: \Lambda_{1}\to \Lambda_{2}$
such that the diagram
$$\begin{CD}
0@>>> W @>>g_{1}> \Lambda_{1} @>>f_{1}>V @>>>0\\
@. @V1_{W}VV @VhVV @VV1_{V}V\\
0@>>> W @>>g_{2}> \Lambda_{2} @>>f_{2}>V @>>>0,
\end{CD}$$
is commutative.}
\end{defn}

The notion of square-zero extension of $V$ by $W$ is an analogue of the notion of 
square-zero extension of an associative algebra by a bimodule. (see, 
for example, Section 9.3 of \cite{W}).

We have the following result for the second cohomology 
$H^{2}_{\frac{1}{2}}(V, W)$ of $V$ with coefficients in $W$:

\begin{thm}\label{h2-ext}
Let $V$ be a grading-restricted vertex algebra and $W$ 
a $V$-module. Then the set of the equivalence classes of 
square-zero 
extensions of $V$ by $W$ corresponds bijectively to $H^{2}_{\frac{1}{2}}(V, W)$.
\end{thm}

\begin{defn}
{\rm Let $t$ be a complex variable. 
A {\it family of grading-restricted vertex algebras up to the first order in $t$} 
is a $\Z$-graded vector space $V$, 
a family $Y_{t}: V\otimes V \to V((x))$ for $t\in \C$ of linear maps 
of the form $Y_{t}=Y_{0}+t\Psi$ where  $Y_{0}$ and
$\Psi$ are linear maps from $V\otimes V$ to $V((x))$ independent of $t$,
and an element $\one\in V$ satisfying all the axioms for grading-restricted vertex algebras
up to  the first order in $t$. More precisely, the triple
$(V, Y_{t}, \one)$ satisfies the grading restriction condition, lower-truncation 
condition for vertex operators, $L(0)$-bracket formula and
the following conditions:
\begin{enumerate}

\item {\it Identity property up to the first order in $t$}: 
$Y_{t}(\one, x)=1_{V}+O(t^{2})$.

\item {\it Creation property  up to the first order in $t$}: For $u\in V$, 
$Y_{t}(u, x)\one\in V[[x]]$ and $\lim_{x\to 0}Y_{t}(u, x)\one=v+O(t^{2})$.

\item {\it Duality up to the first order in $t$}: For $v_{1}, v_{2}, v_{3}\in V$ and 
$v'\in V'$, the coefficients of $t^{0}$ and $t^{1}$ terms of
\begin{eqnarray*}
&\langle v', Y_{t}(v_{1}, z_{1})Y_{t}(v_{2}, z_{2})v_{3}\rangle&\\
&\langle v', Y_{t}(v_{2}, z_{2})Y_{t}(v_{1}, z_{1})v_{3}\rangle&\\
&\langle v', Y_{t}(Y_{t}(v_{1}, z_{1}-z_{2})v_{2}, z_{2})v_{3}\rangle&
\end{eqnarray*}
are absolutely convergent
in the regions $|z_{1}|>|z_{2}|>0$, $|z_{2}|>|z_{1}|>0$ and
$|z_{2}|>|z_{1}-z_{2}|>0$, respectively, to common rational functions 
in $z_{1}$ and $z_{2}$ with the only possible poles at $z_{1}, z_{2}=0$ and 
$z_{1}=z_{2}$.

\end{enumerate}
}
\end{defn}

\begin{defn}
{\rm Let $(V, Y_{V}, \one)$ be a grading-restricted vertex algebra. 
A {\it first order deformation of $V$} is a family 
$Y_{t}: V\otimes V\to V((x))$ for $t\in \C$  of linear maps 
of the form $Y_{t}=Y_{V}+t\Psi$ where 
\begin{eqnarray*}
\Psi: V\otimes V&\to& V((x))\nn
v_{1}\otimes v_{2}&\to &\Psi(v_{1}, x)v_{2}
\end{eqnarray*}
is a linear map
such that 
$(V, Y_{t}, \one)$ for $t\in \C$ is a family of grading-restricted 
vertex algebras up to the first order in $t$.
Two first order deformations $Y_{t}^{(1)}$ 
and $Y_{t}^{(2)}$, $t\in \C$,   
of $(V, Y_{V}, \one)$ are {\it equivalent} if there exists a family
$f_{t}: V \to V$, $t\in \C$ 
of linear maps of the form $f_{t}=1_{V}+tg$
where $g: V\to V$ is a linear map preserving the gradings of $V$ such that 
\begin{equation}\label{equiv-def}
f_{t}(Y^{(1)}_{t}(v_{1}, x)v_{2})-Y^{(2)}_{t}(f_{t}(v_{1}), x)f_{t}(v_{2})
\in t^{2}V((x))
\end{equation}
for $v_{1}, v_{2}\in V$.}
\end{defn}

We have:

\begin{thm}\label{deform-ext}
The set of equivalence classes of first order deformations of 
a grading-restricted vertex algebra is in bijection with 
the set of equivalence classes of square-zero extensions 
of $V$ by $V$.
\end{thm}

From Theorems \ref{h2-ext} and \ref{deform-ext}, we obtain immediately
the following result for the second cohomology 
$H^{2}_{\frac{1}{2}}(V, V)$ of $V$ with coefficients in $V$:

\begin{thm}\label{deform}
Let $V$ be a grading-restricted vertex algebra. 
Then the set of the equivalence classes of first order deformations 
of  $V$  correspond bijectively to $H_{\frac{1}{2}}^{2}(V, V)$.
\end{thm}

We prove Theorems \ref{1st-coh}, \ref{h2-ext} and \ref{deform-ext}
in Sections 2, 3 and 4, respectively.

\paragraph{Acknowledgments}
The author is grateful for partial support {}from NSF
grant PHY-0901237.

\renewcommand{\theequation}{\thesection.\arabic{equation}}
\renewcommand{\thethm}{\thesection.\arabic{thm}}
\setcounter{equation}{0}
\setcounter{thm}{0}

\section{First cohomologies and spaces of derivations}

We prove  Theorem \ref{1st-coh} in the present section.
First, we need the following:

\begin{lemma}\label{f-one}
Let $f: V\to W$ be a derivation. Then $f(\one)=0$.
\end{lemma}
\pf
By definition,
\begin{eqnarray*}
f(\one)&=&f(Y_{V}(\one, z)\one)\nn
&=&\lim_{z\to 0}f(Y_{V}(\one, z)\one)\nn
&=&\lim_{z\to 0}e^{zL(-1)}Y_{W}(\one, -z)f(\one)+
\lim_{z\to 0}Y_{W}(\one, z)f(\one)\nn
&=&2f(\one).
\end{eqnarray*}
So $f(\one)=0$.
\epfv

Let $\Phi: V\to \widetilde{W}_{z_{1}}$ be an element of $C_{m}^{1}(V, W)$ 
satisfying $\delta_{m}^{1} \Phi=0$. Since $\Phi$ satisfies the $L(0)$-conjugation
property, for $v\in V_{(n)}$ and $z\in \C^{\times}$,
\begin{eqnarray*}
z^{L(0)}(\Phi(v))(0)&=&(\Phi(z^{L(0)}v))(0)\nn
&=&z^{n}(\Phi(v))(0).
\end{eqnarray*}
Thus $(\Phi(v))(0)\in W_{(n)}$. So $(\Phi(v))(0)$ is a grading-preserving
linear map from 
$V$ to $W$. 

Since $\delta_{m}^{1} \Phi=0$,
\begin{eqnarray*}
\lefteqn{R(\langle w', 
Y_{W}(v_{1}, z_{1})(\Phi(v_{2}))(z_{2})\rangle)
-R(\langle w', (\Phi(Y_{V}(v_{1}, z_{1}-z_{2})v_{2}))(z_{2})\rangle)}\nn
&&
\quad +R(\langle w', Y_{W}(v_{2}, z_{2})(\Phi(v_{1}))(z_{1})\rangle)\nn
&&=0\quad\quad\quad\quad\quad\quad\quad\quad \quad \quad \quad\quad
\quad\quad\quad\quad\quad\quad\quad \quad \quad \quad
\quad\quad\quad
\end{eqnarray*}
for $v_{1}, v_{2}\in V$ and $w'\in W'$. By $L(-1)$-derivative property for $\Phi$ and 
the vertex operator map $Y_{W}$, 
$$R(\langle w', Y_{W}(v_{2}, z_{2})(\Phi(v_{1}))(z_{1})\rangle)
=R(\langle w', e^{z_{1}L(-1)}Y_{W}(v_{2}, -z_{1}+z_{2})(\Phi(v_{1}))(0)\rangle).$$
Thus we have
\begin{eqnarray*}
\lefteqn{R(\langle w', 
Y_{W}(v_{1}, z_{1})(\Phi(v_{2}))(z_{2})\rangle)
-R(\langle w', (\Phi(Y_{V}(v_{1}, z_{1}-z_{2})v_{2}))(z_{2})\rangle)}\nn
&&\quad +R(\langle w', e^{z_{1}L(-1)}Y_{W}(v_{2}, -z_{1}+z_{2})(\Phi(v_{1}))(0)\rangle)\nn
&&=0.\quad\quad\quad\quad\quad\quad\quad\quad \quad \quad \quad\quad
\quad\quad\quad\quad\quad\quad\quad \quad \quad \quad
\quad\quad\quad
\end{eqnarray*}
Let $z_{2}=0$. We obtain
\begin{eqnarray*}
\lefteqn{R(\langle w', 
Y_{W}(v_{1}, z_{1})(\Phi(v_{2}))(0)\rangle)
-R(\langle w', (\Phi(Y_{V}(v_{1}, z_{1})v_{2}))(0)\rangle)}\nn
&&\quad +R(\langle w', e^{z_{1}L(-1)}Y_{W}(v_{2}, -z_{1})(\Phi(v_{1}))(0)\rangle)\nn
&&=0.\quad\quad\quad\quad\quad\quad\quad\quad \quad \quad \quad\quad
\quad\quad\quad\quad\quad\quad\quad \quad \quad \quad
\quad\quad\quad
\end{eqnarray*}
Since $w'$ is arbitrary, we obtain
\begin{eqnarray*}
\lefteqn{(\Phi(Y_{V}(v_{1}, z_{1})v_{2}))(0)}\nn
&&=e^{z_{1}L(-1)} Y_{W}(v_{2}, -z_{1})(\Phi(v_{1}))(0)
+Y_{W}(v_{1}, z_{1})(\Phi(v_{2}))(0)\nn
&&=Y_{WV}^{W}((\Phi(v_{1}))(0), z_{1})(\Phi(v_{2}))(0)
+Y_{W}(v_{1}, z_{1})(\Phi(v_{2}))(0)
\end{eqnarray*}
for $v_{1}, v_{2}\in V$. This means that $(\Phi(\cdot))(0): V\to W$ 
is a derivation from $V$ to $W$. Note that 
$\delta_{m}^{0}(C_{m}^{0}(V, W))=0$.
So we obtain a linear map from 
$H^{1}(V, W)$ to the space of derivations from $V$ to $W$. 

Conversely, given any derivation $f$ from $V$ to $W$, 
let $\Phi_{f}: V\to \widetilde{W}_{z_{1}}$ be given by 
$$(\Phi_{f}(v))(z_{1})=f(Y_{V}(v, z_{1})\one)=Y_{WV}^{W}(f(v), z_{1})\one$$
for $v\in V$, where we have used Lemma \ref{f-one}. By Theorem 5.6.2 in \cite{FHL}, 
the map from $V$ to $\widetilde{W}_{z_{1}}$ given by 
$v\mapsto Y_{WV}^{W}((\Phi(v))(0), z_{1})\one$ is composable with 
m vertex operators for any $m\in \N$. 
Thus $\Phi_{f}\in C_{m}^{1}(V, W)$ for any $m\in \N$. 
For $v_{1}, v_{2}\in V$ and 
$w'\in W'$, 
\begin{eqnarray}\label{d-phi-f}
\lefteqn{((\delta_{m}^{1}\Phi_{f})(v_{1}\otimes v_{2}))(z_{1}, z_{2})}\nn
&&=R(\langle w', 
Y_{W}(v_{1}, z_{1})Y_{WV}^{W}(f(v_{2}), z_{2})\one\rangle)
\nn
&&\quad -R(\langle w', Y_{WV}^{W}(f(Y_{V}(v_{1}, z_{1}-z_{2})v_{2}), z_{2})\one)
\rangle)\nn
&&
\quad +R(\langle w', Y_{W}(v_{2}, z_{2})Y_{WV}^{W}(f(v_{1}), z_{1})\one\rangle)\nn
&&=R(\langle w', 
Y_{W}(v_{1}, z_{1})Y_{WV}^{W}(f(v_{2}), z_{2})\one\rangle)
\nn
&&\quad -R(\langle w', Y_{WV}^{W}(Y_{WV}^{W}(f(v_{1}), z_{1}-z_{2})v_{2}), z_{2})\one)
\rangle)\nn
&&\quad -R(\langle w', Y_{WV}^{W}(Y_{W}(v_{1}, z_{1}-z_{2})f(v_{2}), z_{2})\one)
\rangle)\nn
&&
\quad +R(\langle w', Y_{W}(v_{2}, z_{2})Y_{WV}^{W}(f(v_{1}), z_{1})\one\rangle)\nn
&&=R(\langle w', 
Y_{W}(v_{1}, z_{1})Y_{WV}^{W}(f(v_{2}), z_{2})\one\rangle)
\nn
&&\quad -R(\langle w', e^{z_{2}L_{W}(-1)}
Y_{WV}^{W}(f(v_{1}), z_{1}-z_{2})v_{2}
\rangle)\nn
&&\quad -R(\langle w', e^{z_{2}L_{W}(-1)}Y_{W}(v_{1}, z_{1}-z_{2})f(v_{2})
\rangle)\nn
&&
\quad +R(\langle w', Y_{W}(v_{2}, z_{2})Y_{WV}^{W}(f(v_{1}), z_{1})\one\rangle)\nn
&&=R(\langle w', 
Y_{W}(v_{1}, z_{1})Y_{WV}^{W}(f(v_{2}), z_{2})\one\rangle)
\nn
&&\quad -R(\langle w', 
Y_{WV}^{W}(f(v_{1}), z_{1})e^{z_{2}L_{V}(-1)}v_{2}
\rangle)\nn
&&\quad -R(\langle w', Y_{W}(v_{1}, z_{1})e^{z_{2}L_{W}(-1)}f(v_{2})
\rangle)\nn
&&
\quad +R(\langle w', Y_{W}(v_{2}, z_{2})Y_{WV}^{W}(f(v_{1}), z_{1})\one\rangle)\nn
&&=R(\langle w', 
Y_{W}(v_{1}, z_{1})Y_{WV}^{W}(f(v_{2}), z_{2})\one\rangle)
\nn
&&\quad -R(\langle w', 
Y_{WV}^{W}(f(v_{1}), z_{1})Y_{W}(v_{2}, z_{2})\one
\rangle)\nn
&&\quad -R(\langle w', Y_{W}(v_{1}, z_{1})Y_{WV}^{W}(f(v_{2}), z_{2})\one
\rangle)\nn
&&
\quad +R(\langle w', Y_{W}(v_{2}, z_{2})Y_{WV}^{W}(f(v_{1}), z_{1})\one\rangle)\nn
&&=-R(\langle w', Y_{WV}^{W}(f(v_{1}), z_{1})Y_{V}(v_{2}, z_{2})\one\rangle)\nn
&&\quad
+R(\langle w', Y_{W}(v_{2}, z_{2})Y_{WV}^{W}(f(v_{1}), z_{1})\one\rangle).\nn
\end{eqnarray}
From Theorem 5.6.2 in \cite{FHL}, we know that the right-hand side of 
(\ref{d-phi-f}) is $0$. 
So we obtain a linear map $f\mapsto \Phi_{f}$ from the space 
$\mbox{\rm Der}\;(V, W)$ to
$H_{m}^{1}(V, W)=C_{m}^{1}(V, W)$. 

Clearly these two maps are inverse to each other and thus $\mbox{\rm Der}\;(V, W)$
and $H_{m}^{1}(V, W)$ are isomorphic.
\epfv

\renewcommand{\theequation}{\thesection.\arabic{equation}}
\renewcommand{\thethm}{\thesection.\arabic{thm}}
\setcounter{equation}{0}
\setcounter{thm}{0}

\section{Second cohomologies and square-zero extensions}

In this section, we prove Theorem \ref{h2-ext}.

Let $(\Lambda, f, g)$  be a square-zero extension of $V$ by $W$. 
Then there is an injective linear map  $\Gamma: V\to \Lambda$
such that the linear map $h: V\oplus W\to \Lambda$ given by
$h(v, w)=\Gamma(v)+g(w)$ is a linear isomorphism. By definition,
the restriction of $h$ to $W$ is the isomorphism $g$
from $W$ to $\ker f$. Then the grading-restricted vertex algebra structure
and the $V$-module structure on $\Lambda$ 
give a grading-restricted vertex algebra structure and a
$V$-module structure on $V\oplus W$ such that the embedding 
$i_{2}: W\to V\oplus W$  and 
the projection $p_{1}: V\oplus W\to V$
are  homomorphisms of grading-restricted vertex algebras. Moreover,
$\ker p_{1}$ is a square-zero ideal of $V\oplus W$, $i_{2}$ is an injective homomorphism
such that $i_{2}(W)=\ker p_{1}$ and
the diagram 
\begin{equation}
\begin{CD}
0@>>> W @>i_{2}>> V\oplus W @>p_{1}>>V @>>>0\\
@. @V1_{W}VV @VhVV @VV1_{V}V\\
0@>>> W @>>g> \Lambda @>>f>V @>>>0
\end{CD}
\end{equation}
of $V$-modules 
is commutative. 
So we obtain a  square-zero extension $(V\oplus W, p_{1}, i_{2})$ 
equivalent to $(\Lambda, f, g)$. We need only consider square-zero extension
of $V$ by $W$ of the particular form $(V\oplus W, p_{1}, i_{2})$. 
Note that the difference between 
two such square-zero extensions are in the vertex operator maps. So 
we use $(V\oplus W, Y_{V\oplus W}, p_{1}, i_{2})$ 
to denote such a 
square-zero extension.

We now write down the vertex operator map for $V\oplus W$ explicitly.
Since $(V\oplus W, Y_{V\oplus W}, p_{1}, i_{2})$ is a square-zero extension of $V$,
there exists $\Psi(u, x)v\in W((x))$ for $u, v\in V$ such that 
\begin{eqnarray*}
Y_{V\oplus W}((v_{1}, 0), x)(v_{2}, 0)&=&(Y_{V}(v_{1}, x)v_{2}, \Psi(v_{1}, x)v_{2}),\\
Y_{V\oplus W}((v_{1}, 0), x)(0, w)&=&(0, Y_{V}(v_{1}, x)w_{2}),\\
Y_{V\oplus W}((0, w_{1}), x)(v_{2}, 0)&=&(0, Y^{W}_{WV}(w, x)v_{2}),\\
Y_{V\oplus W}((0, w_{1}), x)(0, w_{2})&=&0
\end{eqnarray*}
for $v_{1}, v_{2}\in V$ and $w_{1}, w_{2}\in W$.
Thus we have 
\begin{eqnarray}\label{Y_V+W}
\lefteqn{Y_{V\oplus W}((v_{1}, w_{1}), x)(v_{2}, w_{2})}\nn
&&=(Y_{V}(v_{1}, x)v_{2}, 
Y_{W}(v_{1}, x)w_{2}+Y_{WV}^{W}(w_{1}, x)v_{2}+\Psi(v_{1}, x)v_{2})
\end{eqnarray}
for $v_{1}, v_{2}\in V$ and $w_{1}, w_{2}\in W$. 

The vacuum of $V\oplus W$ is $(\one, 0)$. 
Since
\begin{eqnarray*}
Y_{V\oplus W}((v, w), x)(\one, 0)&=&e^{xL_{V\oplus W}(-1)}(v, w)\nn
&=&(e^{xL_{V}(-1)}v, e^{xL_{W}(-1)}w)\nn
&=&(Y_{V}(v, x)\one, Y_{WV}^{W}(w, x)\one)
\end{eqnarray*}
for $v\in W$ and $w\in W$, we have
\begin{equation}\label{psi-one}
\Psi(v, x)\one=0
\end{equation}
for $v\in V$.

We identify $(V\oplus W)'$ with 
$V'\oplus W'$. 
For $v_{1}, v_{2}\in V$ and $w'\in W'$, 
\begin{eqnarray*}
\lefteqn{\langle (0, w'), Y_{V\oplus W}((v_{1}, 0), z_{1})Y_{V\oplus W}((v_{2}, 0), z_{2})
(\one, 0)\rangle}\nn
&&=\langle w', \Psi(v_{1}, z_{1})Y_{V}(v_{2}, z_{2})\one
+Y_{W}(v_{1}, z_{1})\Psi(v_{2}, z_{2})\one\rangle\nn
&&=\langle w', \Psi(v_{1}, z_{1})Y_{V}(v_{2}, z_{2})\one\rangle,\\
\lefteqn{\langle (0, w'), Y_{V\oplus W}((v_{2}, 0), z_{2})Y_{V\oplus W}((v_{1}, 0), z_{1})
(\one, 0)\rangle}\nn
&&=\langle w', \Psi(v_{2}, z_{2})Y_{V}(v_{1}, z_{1})\one
+Y_{W}(v_{2}, z_{2})\Psi(v_{1}, z_{1})\one\rangle\nn
&&=\langle w', \Psi(v_{2}, z_{2})Y_{V}(v_{1}, z_{1})\one\rangle,\\
\lefteqn{\langle (0, w'), Y_{V\oplus W}(Y_{V\oplus W}((v_{1}, 0), z_{1}-z_{2})
(v_{2}, 0), z_{2})
(\one, 0)\rangle}\nn
&&=\langle w', Y_{WV}^{W}(\Psi(v_{1}, z_{1}-z_{2})
v_{2}, z_{2})\one+\Psi(Y_{V}(v_{1}, z_{1}-z_{2})
v_{2}, z_{2})\one\rangle\nn
&&=\langle w', Y_{WV}^{W}(\Psi(v_{1}, z_{1}-z_{2})
v_{2}, z_{2})\one\rangle
\end{eqnarray*}
are absolutely convergent in the region $|z_{1}|>|z_{2}|>0$, $|z_{2}|>|z_{1}|>0$, 
$|z_{2}|>|z_{1}-z_{2}|>0$, respectively, to one
rational function in $z_{1}$ and $z_{2}$
with the only possible poles at $z_{1}, z_{2}=0$ and $z_{1}=z_{2}$. 
Using our notation in \cite{H}, we denote this rational function by
$$R(\langle w', \Psi(v_{1}, z_{1})Y_{V}(v_{2}, z_{2})\one\rangle)$$
or 
$$R(\langle w', \Psi(v_{2}, z_{2})Y_{V}(v_{1}, z_{1})\one\rangle)$$
or 
$$R(\langle w', Y_{WV}^{W}(\Psi(v_{1}, z_{1}-z_{2})
v_{2}, z_{2})\one\rangle).$$
Then we obtain an element, denoted by
$$E(\Psi(v_{1}, z_{1})Y_{V}(v_{2}, z_{2})\one)$$
or 
$$E(\Psi(v_{2}, z_{2})Y_{V}(v_{1}, z_{1})\one)$$
or 
$$E(Y_{WV}^{W}(\Psi(v_{1}, z_{1}-z_{2})
v_{2}, z_{2})\one),$$
of $\widetilde{W}_{z_{1}, z_{2}}$
given by 
$$\langle w', E(\Psi(v_{1}, z_{1})Y_{V}(v_{2}, z_{2})\one)\rangle
=R(\langle w', \Psi(v_{1}, z_{1})Y_{V}(v_{2}, z_{2})\one\rangle)$$
or 
$$\langle w', E(\Psi(v_{2}, z_{2})Y_{V}(v_{1}, z_{1})\one)\rangle
=R(\langle w', \Psi(v_{2}, z_{2})Y_{V}(v_{1}, z_{1})\one\rangle)$$
or 
$$\langle w', E(Y_{WV}^{W}(\Psi(v_{1}, z_{1}-z_{2})
v_{2}, z_{2})\one)\rangle
=R(\langle w', Y_{WV}^{W}(\Psi(v_{1}, z_{1}-z_{2})
v_{2}, z_{2})\one\rangle).$$
By definition, we have 
\begin{eqnarray*}
E(\Psi(v_{1}, z_{1})Y_{V}(v_{2}, z_{2})\one)
&=&E(\Psi(v_{2}, z_{2})Y_{V}(v_{1}, z_{1})\one)\nn
&=&E(Y_{WV}^{W}(\Psi(v_{1}, z_{1}-z_{2})
v_{2}, z_{2})\one)
\end{eqnarray*}
for $v_{1}, v_{2}\in V$.

Let 
$$\Phi: V\otimes V\to \widetilde{W}_{z_{1}, z_{2}}$$
be the linear map given by
\begin{eqnarray}\label{phi}
(\Phi(v_{1}\otimes v_{2}))(z_{1}, z_{2})
&=&E(\Psi(v_{1}, z_{1})Y_{V}(v_{2}, z_{2})\one)\nn
&=&E(\Psi(v_{2}, z_{2})Y_{V}(v_{1}, z_{1})\one)\nn
&=&E(Y_{WV}^{W}(\Psi(v_{1}, z_{1}-z_{2})
v_{2}, z_{2})\one)
\end{eqnarray}
for $v_{1}, v_{2}\in V$ and $(z_{1}, z_{2})\in F_{2}\C$. 
We first prove that $\Phi\in \widehat{C}_{\frac{1}{2}}^{2}(V, W)$.

By the $L(-1)$-derivative property and the $L(0)$-bracket formula 
for $V\oplus W$, we have 
\begin{eqnarray}
\frac{d}{dx}Y_{V\oplus W}((v, 0), x)
&=&Y_{V\oplus W}((L_{V}(-1)v, 0), x)\label{l-1l0.1}\\
&=&[L_{V+W}(-1), Y_{V\oplus W}((v, 0), x)],\label{l-1l0.2}\\
{[L_{V+W}(0), Y_{V\oplus W}((v, 0), x)]}
&=&Y_{V\oplus W}((L_{V}(0)v, 0), x)
+x\frac{d}{dx}Y_{V\oplus W}((v, 0), x)\label{l-1l0.3}\nn
\end{eqnarray}
for $v\in V$. By (\ref{l-1l0.1}), (\ref{l-1l0.2}),
(\ref{l-1l0.3}), (\ref{Y_V+W}) and the 
$L(-1)$-derivative property and the $L(0)$-bracket formula 
for $V$, we obtain
\begin{eqnarray}
\frac{d}{dx}\Psi(v, x)
&=&\Psi(L_{V}(-1)v, x)\label{psi-l-1l0.1}\\
&=&L_{W}(-1)\Psi(v, x)-\Psi(v, x)L_{V}(-1),\label{psi-l-1l0.2}\nn
&&\\
L_{W}(0)\Psi(v, x)
-\Psi(v, x)L_{V}(0)
&=&\Psi(L_{V}(0)v, x)
+x\frac{d}{dx}\Psi(v, x)\label{psi-l-1l0.3}
\end{eqnarray}
for $v\in V$. From (\ref{psi-l-1l0.3}), we obtain
\begin{equation}\label{conjugation}
z^{L_{W}(0)}\Psi(v, x)=\Psi(z^{L_{V}(0)}v, zx)z^{L_{V}(0)}
\end{equation}
for $v\in V$.

For $v_{1}, v_{2}\in V$ and $w'\in W'$, 
by (\ref{psi-l-1l0.1}) and the $L(-1)$-derivative property for $V$,
we obtain
\begin{eqnarray}\label{partial-z1}
\lefteqn{\frac{\partial}{\partial z_{1}}\langle w', 
(\Phi(v_{1}\otimes v_{2}))(z_{1}, z_{2})
\rangle}\nn
&&=\frac{\partial}{\partial z_{1}}\langle w', 
E(\Psi(v_{1}, z_{1})Y_{V}(v_{2}, z_{2})\one)
\rangle\nn
&&=\frac{\partial}{\partial z_{1}}R(\langle w', 
\Psi(v_{1}, z_{1})Y_{V}(v_{2}, z_{2})\one
\rangle)\nn
&&=R\left(\left\langle w', 
\frac{\partial}{\partial z_{1}}\Psi(v_{1}, z_{1})Y_{V}(v_{2}, z_{2})
\one\right\rangle\right)\nn
&&=R(\langle w', 
\Psi(L_{V}(-1)v_{1}, z_{1})Y_{V}(v_{2}, z_{2})\one\rangle)\nn
&&=\langle w', 
E(\Psi(L_{V}(-1)v_{1}, z_{1})Y_{V}(v_{2}, z_{2})\one)\rangle\nn
&&=\langle w', (\Phi(L_{V}(-1)v_{1}\otimes v_{2}))(z_{1}, z_{2})
\rangle
\end{eqnarray}
and
\begin{eqnarray}\label{partial-z2}
\lefteqn{\frac{\partial}{\partial z_{2}}\langle w', 
(\Phi(v_{1}\otimes v_{2}))(z_{1}, z_{2})
\rangle}\nn
&&=\frac{\partial}{\partial z_{2}}\langle w', 
E(\Psi(v_{1}, z_{1})Y_{V}(v_{2}, z_{2})\one)
\rangle\nn
&&=\frac{\partial}{\partial z_{2}}R(\langle w', 
\Psi(v_{1}, z_{1})Y_{V}(v_{2}, z_{2})\one
\rangle)\nn
&&=R\left(\left\langle w', 
\Psi(v_{1}, z_{1})\frac{\partial}{\partial z_{2}}Y_{V}(v_{2}, z_{2})
\one\right\rangle\right)\nn
&&=R(\langle w', 
\Psi(v_{1}, z_{1})Y_{V}(L_{V}(-1)v_{2}, z_{2})\one\rangle)\nn
&&=\langle w', 
E(\Psi(v_{1}, z_{1})Y_{V}(L_{V}(-1)v_{2}, z_{2})\one)\rangle\nn
&&=\langle w', (\Phi(v_{1}\otimes L_{V}(-1)v_{2}))(z_{1}, z_{2})
\rangle.
\end{eqnarray}
Using (\ref{psi-l-1l0.1}), (\ref{psi-l-1l0.2})
and the $L(-1)$-derivative property for $V$, we obtain
\begin{eqnarray}\label{partial-z1z2}
\lefteqn{\left(\frac{\partial}{\partial z_{2}}+
\frac{\partial}{\partial z_{2}}\right)
\langle w', (\Phi(v_{1}\otimes v_{2}))(z_{1}, z_{2})
\rangle}\nn
&&=\left(\frac{\partial}{\partial z_{2}}+
\frac{\partial}{\partial z_{2}}\right)\langle w', 
E(\Psi(v_{1}, z_{1})Y_{V}(v_{2}, z_{2})\one)
\rangle\nn
&&=\left(\frac{\partial}{\partial z_{2}}+
\frac{\partial}{\partial z_{2}}\right)R(\langle w', 
\Psi(v_{1}, z_{1})Y_{V}(v_{2}, z_{2})\one)
\rangle)\nn
&&=R\left(\left\langle w', 
\frac{\partial}{\partial z_{1}}\Psi(v_{1}, z_{1})Y_{V}(v_{2}, z_{2})
\one\right\rangle\right)\nn
&&\quad+R\left(\left\langle w', 
\Psi(v_{1}, z_{1})\frac{\partial}{\partial z_{2}}Y_{V}(v_{2}, z_{2})
\one\right\rangle\right)\nn
&&=R(\langle w', 
\Psi(L_{V}(-1)v_{1}, z_{1})Y_{V}(v_{2}, z_{2})\one\rangle)\nn
&&\quad +R(\langle w', 
\Psi(v_{1}, z_{1})Y_{V}(L_{V}(-1)v_{2}, z_{2})\one\rangle)\nn
&&=R(\langle w', 
L_{W}(-1)\Psi(v_{1}, z_{1})Y_{V}(v_{2}, z_{2})\one\rangle)\nn
&&=R(\langle L_{W'}(1)w', 
\Psi(v_{1}, z_{1})Y_{V}(v_{2}, z_{2})\one\rangle)\nn
&&=\langle L_{W'}(1)w', 
E(\Psi(v_{1}, z_{1})Y_{V}(v_{2}, z_{2})\one)\rangle\nn
&&=\langle w', L_{W}(-1)
E(\Psi(v_{1}, z_{1})Y_{V}(v_{2}, z_{2})\one)\rangle\nn
&&=\langle w', L_{W}(-1)(\Phi(v_{1}\otimes v_{2}))(z_{1}, z_{2})
\rangle
\end{eqnarray}
for $v_{1}, v_{2}\in V$ and $w'\in W'$, 
From (\ref{partial-z1}), (\ref{partial-z2}) and (\ref{partial-z1z2}),
we see that $\Phi$ satisfies the $L(-1)$-derivative property. 

Also for $v_{1}, v_{2}\in V$ and $w'\in W'$, by (\ref{conjugation}) and 
the $L(0)$-bracket formula for $V$, we have
\begin{eqnarray*}
\lefteqn{\langle w', z^{L_{W}(0)}
(\Phi(v_{1}\otimes v_{2}))(z_{1}, z_{2})
\rangle}\nn
&&=\langle w', z^{L_{W}(0)}E(\Psi(v_{1}, z_{1})Y_{V}(v_{2}, z_{2})\one)
\rangle\nn
&&=\langle z^{L_{W'}(0)}w', E(\Psi(v_{1}, z_{1})Y_{V}(v_{2}, z_{2})\one)
\rangle\nn
&&=R(\langle z^{L_{W'}(0)}w', \Psi(v_{1}, z_{1})Y_{V}(v_{2}, z_{2})\one
\rangle)\nn
&&=R(\langle w', z^{L_{W}(0)}\Psi(v_{1}, z_{1})Y_{V}(v_{2}, z_{2})\one
\rangle)\nn
&&=R(\langle w', \Psi(z^{L_{V}(0)}v_{1}, zz_{1})Y_{V}(z^{L_{V}(0)}v_{2}, zz_{2})\one
\rangle)\nn
&&=\langle w', E(\Psi(z^{L_{V}(0)}v_{1}, zz_{1})Y_{V}(z^{L_{V}(0)}v_{2}, zz_{2})\one)
\rangle\nn
&&=\langle w', 
(\Phi(z^{L_{V}(0)}v_{1}\otimes z^{L_{V}(0)}v_{2}))(zz_{1}, zz_{2})
\rangle,
\end{eqnarray*}
that is, $\Phi$ satisfies the $L(0)$-conjugation property.

Since $V\oplus W$ is a grading-restricted vertex algebra, 
for $v_{1}, v_{2}, v_{3}\in V$ and $w'\in W'$, the series 
$$\langle (0, w'), Y_{V\oplus W}((v_{1}, 0), z_{1})
Y_{V\oplus W}((v_{2}, 0), z_{2})Y_{V\oplus W}((v_{3}, 0), z_{3})(\one, 0)\rangle$$
and
$$\langle (0, w'), Y_{V\oplus W}(Y_{V\oplus W}((v_{1}, 0), z_{1}-z_{2})
(v_{2}, 0), 
z_{2})Y_{V\oplus W}((v_{3}, 0), z_{3})(\one, 0)\rangle$$
are absolutely convergent in the regions given by $|z_{1}|>|z_{2}|> 
|z_{3}|>0$ and
by $|z_{2}|>|z_{1}-z_{2}|, |z_{3}|>0$
and $|z_{2}-z_{3}|>|z_{1}-z_{2}|$, respectively, 
to a same rational function
with the only possible poles at $z_{1}=z_{2}$,
$z_{1}=z_{3}$, $z_{2}=z_{3}$. 
But by (\ref{Y_V+W}) and (\ref{psi-one}), these series are equal to 
$$\langle w', \Psi(v_{1}, z_{1})
Y_{V}(v_{2}, z_{2})Y_{V}(v_{3}, z_{3})\one\rangle+\langle w', Y_{W}(v_{1}, z_{1})
\Psi(v_{2}, z_{2})Y_{V}(v_{3}, z_{3})\one\rangle$$
and 
\begin{eqnarray*}
\lefteqn{\langle w', \Psi(Y_{V}(v_{1}, z_{1}-z_{2})v_{2}, 
z_{2})Y_{V}(v_{3}, z_{3})\one\rangle}\nn
&&\quad 
+\langle w', Y_{WV}^{W}(\Psi(v_{1}, z_{1}-z_{2})v_{2}, 
z_{2})Y_{V}(v_{3}, z_{3})\one\rangle,
\end{eqnarray*}
respectively, and are absolutely convergent to a same rational function which 
in our convention is equal to 
$$R(\langle w', \Psi(v_{1}, z_{1})
Y_{V}(v_{2}, z_{2})Y_{V}(v_{3}, z_{3})\one\rangle+\langle w', Y_{W}(v_{1}, z_{1})
\Psi(v_{2}, z_{2})Y_{V}(v_{3}, z_{3})\one\rangle)$$
and
\begin{eqnarray*}
\lefteqn{R(\langle w', \Psi(Y_{V}(v_{1}, z_{1}-z_{2})v_{2}, 
z_{2})Y_{V}(v_{3}, z_{3})\one\rangle}\nn
&&\quad \quad\quad\quad\quad\quad 
+\langle w', Y_{WV}^{W}(\Psi(v_{1}, z_{1}-z_{2})v_{2}, 
z_{2})Y_{V}(v_{3}, z_{3})\one\rangle).
\end{eqnarray*}
In particular,  we have proved that 
$\Phi\in \widehat{C}_{\frac{1}{2}}^{2}(V, W)$.

Since by (\ref{phi}),
\begin{eqnarray*}
(\Phi(v_{1}\otimes v_{2}))(z_{1}, z_{2})&=&E(\Psi(v_{1}, z_{1})Y_{V}(v_{2}, z_{2})\one)\nn
&=&E(\Psi(v_{2}, z_{2})Y_{V}(v_{1}, z_{1})\one)\nn
&=&(\Phi(v_{2}\otimes v_{1}))(z_{2}, z_{1})\nn
&=&(\sigma_{12}(\Phi(v_{2}\otimes v_{1})))(z_{1}, z_{2})
\end{eqnarray*}
for $v_{1}, v_{2}\in V$ and $(z_{1}, z_{2})\in F_{2}\C$, that is,
$$\Phi(v_{1}\otimes v_{2})-\sigma_{12}(\Phi(v_{2}\otimes v_{1}))=0$$
for $v_{1}, v_{2}\in V$,
We obtain
\begin{eqnarray*}
\lefteqn{\sum_{\sigma\in J_{2; 1}}(-1)^{|\sigma|}\sigma(\Phi(v_{1}\otimes v_{2}))}\nn
&&=\Phi(v_{1}\otimes v_{2})-\sigma_{12}(\Phi(v_{2}\otimes v_{1}))\nn
&&=0
\end{eqnarray*}
for $v_{1}, v_{2}\in V$. So $\Phi\in C^{2}(V, W)$. 

Next we show that $\delta_{\frac{1}{2}}^{2}(\Phi)=0$. 
For $v_{1}, v_{2}, v_{3}\in V$, $w'\in W'$,
\begin{eqnarray}\label{cocycle}
\lefteqn{\langle w', ((\delta_{\frac{1}{2}}^{2}(\Phi))(v_{1}\otimes v_{2}\otimes v_{3}))
(z_{1}, z_{2}, z_{3})\rangle}\nn
&&=R(\langle w', (E^{(1)}_{W}(v_{1}; \Phi(v_{2}\otimes v_{3})))(z_{1}, z_{2}, z_{3})\rangle
\nn
&&\quad\quad\quad +\langle w', (\Phi(v_{1}\otimes E^{(2)}(v_{2}\otimes v_{3}; \one)))
(z_{1}, z_{2}, z_{3})\rangle)
\nn
&&\quad-R(\langle w', (\Phi(E^{(2)}(v_{1}\otimes v_{2}; \one)\otimes v_{3}))
(z_{1}, z_{2}, z_{3})\rangle\nn
&&\quad\quad\quad +\langle w', (E^{W; (1)}_{WV}(\Phi(v_{1}\otimes v_{2}); v_{3}))
(z_{1}, z_{2}, z_{3})\rangle)\nn
&&=R(\langle w', Y_{W}(v_{1}, z_{1})\Psi(v_{2}, z_{2})Y_{V}(v_{3}, z_{3})\one\rangle
\nn
&&\quad \quad\quad+\langle w', \Psi(v_{1}, z_{1})Y_{V}(v_{2}, z_{2})Y_{V}(v_{3}, z_{3})\one\rangle)
\nn
&&\quad-R(\langle w', \Psi(Y_{V}(v_{1}, z_{1}-z_{2})v_{2},
z_{2})Y_{V}(v_{3}, z_{3})\one\rangle\nn
&&\quad \quad\quad +\langle w', Y_{WV}^{W} 
(\Psi(v_{1}, z_{1}-z_{2})v_{2}, z_{2})Y_{V}(v_{3}, z_{3})\one\rangle).
\end{eqnarray}
Since $V\oplus W$ is a grading-restricted vertex algebra, 
we have the associativity property
\begin{eqnarray*}
\lefteqn{R(\langle (0, w'), Y_{V\oplus W}((v_{1}, 0), z_{1})
Y_{V\oplus W}((v_{2}, 0), z_{2})Y_{V\oplus W}((v_{3}, 0), z_{3})(\one, 0)\rangle)}\nn
&&=R(\langle (0, w'), Y_{V\oplus W}(Y_{V\oplus W}((v_{1}, 0), z_{1}-z_{2})(v_{2}, 0), 
z_{2})\cdot\nn
&&\quad\quad\quad\quad\quad\quad\quad\quad\quad\quad\quad\quad
\quad\quad\quad\quad\quad
\cdot  Y_{V\oplus W}((v_{3}, 0), z_{3})(\one, 0)\rangle),
\end{eqnarray*}
which, by (\ref{Y_V+W}) and (\ref{psi-one}), is equivalent to 
\begin{eqnarray*}
\lefteqn{R(\langle w', \Psi(v_{1}, z_{1})
Y_{V}(v_{2}, z_{2})Y_{V}(v_{3}, z_{3})\one\rangle}\nn
&&\quad\quad \quad  +\langle w', Y_{W}(v_{1}, z_{1})
\Psi(v_{2}, z_{2})Y_{V}(v_{3}, z_{3})\one\rangle)\nn
&&=R(\langle w', \Psi(Y_{V}(v_{1}, z_{1}-z_{2})v_{2}, 
z_{2})Y_{V}(v_{3}, z_{3})\one\rangle\nn
&&\quad\quad \quad  +\langle w', Y_{WV}^{W}(\Psi(v_{1}, z_{1}-z_{2})v_{2}, 
z_{2})Y_{V}(v_{3}, z_{3})\one\rangle),
\end{eqnarray*}
as we have noticed above. So the right-hand side of (\ref{cocycle}) is $0$.
Thus $\Phi+\delta_{2}^{1}C_{2}^{1}(V, W)$ is an element of 
$H_{\frac{1}{2}}^{2}(V, W)$.

Conversely, given any element of $H^{2}_{\frac{1}{2}}(V, W)$, let 
$\Phi\in C^{2}_{\frac{1}{2}}(V, W)$ be a representative of this 
element. Then for any $v_{1}, v_{2}\in V$, there exists $N$ such that 
for $w'\in W'$, 
$\langle w', (\Phi(v_{1}\otimes v_{2}))(z, 0)\rangle$
is a rational function of $z$ with the only possible pole at
$z=0$
of order less than or equal to $N$. For $v_{1}, v_{2}\in V$,
let $\Psi(v_{1}, x)v_{2}\in W((x))$ be given by 
$$\langle w', \Psi(v_{1}, x)v_{2}\rangle|_{x=z}
=\langle w', (\Phi(v_{1}\otimes v_{2}))(z, 0)\rangle.$$
for $z\in \C^{\times}$. For $v_{1}, v_{2}\in V$, define
$Y_{V\oplus W}(v_{1}, x)v_{2}$ using (\ref{Y_V+W}). So we obtain a 
vertex operator map $Y_{V\oplus W}$.
Reversing the proof above, we see that 
$V\oplus W$ equipped with the vertex operator map $Y_{V+W}$ and
the vacuum $(\one, 0)$ is a grading-restricted 
vertex algebra and together with the projection $p_{1}: V\oplus W\to V$ and the 
embedding $i_{2}: W\to V\oplus W$, $V\oplus W$ is a square-zero extension of 
$V$ by $W$. 

Next we prove that two elements of $\ker \delta_{\frac{1}{2}}^{2}$ 
obtained this way differ by 
an element of $\delta_{1}C^{1}(V, W)$ if and only if
the corresponding square-zero extensions of $V$ by $W$ are
equivalent.

Let $\Phi_{1}, \Phi_{2}\in \ker \delta_{\frac{1}{2}}^{2}$ 
be two such elements obtained from 
square-zero extensions $(V\oplus W, Y_{V\oplus W}^{(1)}, 
p_{1}, i_{2})$ and 
$(V\oplus W, Y_{V\oplus W}^{(2)}, p_{1}, i_{2})$. 
Assume that $\Phi_{1}=\Phi_{2}+\delta_{1}(\Gamma)$
where $\Gamma\in C^{1}(V, W)$. Since 
\begin{eqnarray*}
\lefteqn{\langle w', ((\delta_{1}(\Gamma))(v_{1}\otimes v_{2}))(z_{1}, z_{2})\rangle}\nn
&&=R(\langle w', Y_{W}(v_{1}, z_{1})(\Gamma(v_{2}))(z_{2})\rangle)\nn
&&\quad 
-R(\langle w', (\Gamma(Y_{V}(v_{1}, z_{1}-z_{2})v_{2}))(z_{2})\rangle)\nn
&&\quad 
+R(\langle w', Y_{W}(v_{2}, z_{2})(\Gamma(v_{1}))(z_{1})\rangle),
\end{eqnarray*}
we have 
\begin{eqnarray}\label{h2-ext-4}
\lefteqn{R(\langle w', \Psi_{1}(v_{1}, z_{1})Y_{V}(v_{2}, z_{2})\one\rangle)}\nn
&&=\langle w', (\Phi_{1}(v_{1}\otimes v_{2}))(z_{1}, z_{2})\rangle\nn
&&=\langle w', (\Phi_{2}(v_{1}\otimes v_{2}))(z_{1}, z_{2})\rangle\nn
&&\quad 
+\langle w', (\delta_{1}(\Gamma))(z_{1}, z_{2})\rangle\nn
&&=R(\langle w', \Psi_{2}(v_{1}, z_{1})Y_{V}(v_{2}, z_{2})\one\rangle\nn
&&\quad 
+R(\langle w', Y_{W}(v_{1}, z_{1})(\Gamma(v_{2}))(z_{2})\rangle)\nn
&&\quad 
-R(\langle w', (\Gamma(Y_{V}(v_{1}, z_{1}-z_{2})v_{2}))(z_{2})\rangle)\nn
&&\quad 
+R(\langle w', Y_{W}(v_{2}, z_{2})(\Gamma(v_{1}))(z_{1})\rangle)\nn
&&=R(\langle w', \Psi_{2}(v_{1}, z_{1})Y_{V}(v_{2}, z_{2})\one\rangle\nn
&&\quad 
+R(\langle w', Y_{W}(v_{1}, z_{1})(\Gamma(v_{2}))(z_{2})\rangle)\nn
&&\quad 
-R(\langle w', (\Gamma(Y_{V}(v_{1}, z_{1}-z_{2})v_{2}))(z_{2})\rangle)\nn
&&\quad 
+R(\langle w', e^{(z_{1}-z_{2})L_{W}(-1)}Y_{W}(v_{2}, -z_{1})
(\Gamma(v_{1}))(z_{2})\rangle).
\end{eqnarray}
Let $z_{2}$ go to zero on both sides of (\ref{h2-ext-4}). We obtain
\begin{eqnarray*}
\langle w', \Psi_{1}(v_{1}, z_{1})v_{2}\rangle
&=&\langle w', \Psi_{2}(v_{1}, z_{1})v_{2}\rangle\nn
&&
+\langle w', Y_{W}(v_{1}, z_{1})(\Gamma(v_{2}))(0)\rangle\nn
&&
-\langle w', (\Gamma(Y_{V}(v_{1}, z_{1})v_{2}))(0)\rangle\nn
&&
+\langle w', e^{z_{1}L_{W}(-1)}Y_{W}(v_{2}, -z_{1})(\Gamma(v_{1}))(0)\rangle\nn
&=&\langle w', \Psi_{2}(v_{1}, z_{1})v_{2}\rangle\nn
&&
+\langle w', Y_{W}(v_{1}, z_{1})(\Gamma(v_{2}))(0)\rangle\nn
&&
-\langle w', (\Gamma(Y_{V}(v_{1}, z_{1})v_{2}))(0)\rangle\nn
&&
+\langle w', Y_{WV}^{W}((\Gamma(v_{1}))(0), z_{1})v_{2}\rangle.
\end{eqnarray*}
Then 
\begin{eqnarray}\label{h2-ext-5}
\Psi_{1}(v_{1}, x)v_{2}
&=&\Psi_{2}(v_{1}, x)v_{2}
+Y_{W}(v_{1}, x)(\Gamma(v_{2}))(0)\nn
&&
-(\Gamma(Y_{V}(v_{1}, x)v_{2}))(0)
+Y_{WV}^{W}((\Gamma(v_{1}))(0), x)v_{2}.
\end{eqnarray}

For $v_{1}, v_{2}\in V$ and $w_{1}, w_{2}\in W$, 
by (\ref{Y_V+W}) and (\ref{h2-ext-5}), we have 
\begin{eqnarray}\label{h2-ext-6}
\lefteqn{Y_{V\oplus W}^{(1)}((v_{1}, w_{1}), x)(v_{2}, w_{2})}\nn
&&=(Y_{V}(v_{1}, x)v_{2}, Y_{W}(v_{1}, x)w_{2}
+Y_{WV}^{W}(w_{1}, x)v_{2}+ \Psi_{1}(v_{1}, x)v_{2})\nn
&&=(Y_{V}(v_{1}, x)v_{2}, Y_{W}(v_{1}, x)w_{2}
+ Y_{WV}^{W}(w_{1}, x)v_{2}
+ \Psi_{2}(v_{1}, x)v_{2})\nn
&&\quad +(Y_{V}(v_{1}, x)v_{2}, Y_{W}(v_{1}, x)(\Gamma(v_{2}))(0))
\nn
&&\quad -(Y_{V}(v_{1}, x)v_{2}, (\Gamma(Y_{V}(v_{1}, x)v_{2}))(0))\nn
&&\quad +(Y_{V}(v_{1}, x)v_{2}, Y_{WV}^{W}((\Gamma(v_{1}))(0), x)v_{2})\nn
&&=Y_{V\oplus W}^{(2)}((v_{1}, w_{1}+(\Gamma(v_{1}))(0)), x)(v_{2}, 
w_{2}+(\Gamma(v_{2}))(0))\nn
&&\quad  -(Y_{V}(v_{1}, x)v_{2}, (\Gamma(Y_{V}(v_{1}, x)v_{2}))(0)).
\end{eqnarray}

We now define 
a linear map $h: V\oplus W\to V\oplus W$ by
$$e(v, w)=(v, w+(\Gamma(v))(0))$$
for $v\in V$ and $w\in W$. Then $h$ is a linear isomorphism
and (\ref{h2-ext-6}) can be rewritten as 
\begin{equation}\label{e-map-prop}
h(Y_{V\oplus W}^{(1)}((v_{1}, w_{1}), x)(v_{2}, w_{2}))
=Y_{V\oplus W}^{(2)}(h(v_{1}, w_{1}), x)h(v_{2}, w_{2}).
\end{equation}
for $v_{1}, v_{2}\in V$ and $w_{1}, w_{2}\in W$. 
Thus $h$ is an isomorphism of grading-restricted vertex algebras from 
$(V\oplus W, Y_{V\oplus W}^{(1)}, (\one, 0))$ to $(V\oplus W, Y_{V\oplus W}^{(2)}, 
(\one, 0))$
such that the diagram 
\begin{equation}\label{equivalence}
\begin{CD}
0@>>> W @>i_{2}>> V\oplus W @>p_{1}>>V @>>>0\\
@. @V1_{W}VV @VhVV @VV1_{V}V\\
0@>>> W @>i_{2}>> V\oplus W @>p_{1}>>V @>>>0\\
\end{CD}
\end{equation}
is commutative. Thus the two square-zero extensions of $V$ by $W$ are 
equivalent. 

Conversely, let $(V\oplus W, Y_{V+W}^{(1)}, 
p_{1}, i_{2})$ and 
$(V\oplus W, Y_{V+W}^{(2)}, p_{1}, i_{2})$ be two equivalent 
square-zero extensions of $V$ by $W$. So there exists an isomorphism
$h: V\oplus W\to V\oplus W$ of grading-restricted vertex algebras
such that (\ref{equivalence})
is commutative. We have the following lemma which is also needed 
in the next section:

\begin{lemma}\label{h-form}
There exists a linear map $g: V\to V$ such that 
$$h(v, w)=(v, w+g(v))$$
for $v\in V$ and $w\in W$.
\end{lemma}
\pf
Let $h(v, w)=(f(v, w), g(v, w))$ for $v\in V$ and $w\in W$. Then 
by (\ref{equivalence}), we have
$f(v, w)=v$ and $g(0, w)=w$. Since $h$ is linear, we have 
$g(v, w)=g(v, 0)+g(0, w)=w+g(v, 0)$. So $h(v, w)=(v, w+g(v, 0))$. 
Taking $g(v)$ to be $g(v, 0)$, we see that the conclusion holds.
\epfv

Let $(\Gamma(v))(z_{1})=e^{z_{1}L_{W}(-1)}g(v)\in \overline{W}$. Then 
$\Gamma: V\to \widetilde{W}_{z_{1}}$ is an element of $C_{2}^{1}(V, W)$. 
By definition, we have $g(v)=(\Gamma(v))(0)$ and 
$h(v, w)=(v, w+(\Gamma(v))(0))$ for $v\in V$ and $w\in W$. 
Let $\Phi_{1}$ and $\Phi_{2}$ be elements of $\ker \delta_{\frac{1}{2}}^{2}$ 
obtained from $(V\oplus W, Y_{V\oplus W}^{(1)}, 
p_{1}, i_{2})$ and 
$(V\oplus W, Y_{V+W}^{(2)}, p_{1}, i_{2})$, respectively, and 
$\Psi_{1}$ and $\Psi_{2}$ the corresponding maps from $V\otimes V$ to $W((x))$. 
Then since $h$ is a homomorphism of grading-restricted vertex algebras,  
(\ref{e-map-prop}) holds 
for $v_{1}, v_{2}\in V$ and $w_{1}, w_{2}\in W$.
Thus the two sides of (\ref{h2-ext-6}) are equal 
for $v_{1}, v_{2}\in V$ and $w_{1}, w_{2}\in W$. So the two expressions in the 
middle of (\ref{h2-ext-6}) are equal for $v_{1}, v_{2}\in V$ and $w_{1}, w_{2}\in W$.
Thus we have (\ref{h2-ext-5}) for $v_{1}, v_{2}\in V$. 
Formula (\ref{h2-ext-5}) implies that the two sides of (\ref{h2-ext-4})
are equal for $v_{1}, v_{2}\in V$. Thus the middle expressions in  (\ref{h2-ext-4}) 
are all equal for $v_{1}, v_{2}\in V$. In particular, we obtain 
$\Phi_{1}=\Phi_{2}+\Gamma$. So $\Phi_{1}$ and $\Phi_{2}$ differ by 
an element of 
$\delta_{1}C^{1}(V, W)$.
\epfv

\renewcommand{\theequation}{\thesection.\arabic{equation}}
\renewcommand{\thethm}{\thesection.\arabic{thm}}
\setcounter{equation}{0}
\setcounter{thm}{0}

\section{Square-zero extensions and first order deformations}

In this section, we prove Theorem \ref{deform-ext}.

Let $Y_{t}: V\otimes V\to V((x))$, $t\in U$, be a first order deformation of $V$.
By definition, there exists 
\begin{eqnarray*}
\Psi: V\otimes V&\to& V((x))\nn
v_{1}\otimes v_{2}&\to &\Psi(v_{1}, x)v_{2}
\end{eqnarray*}
such that 
$$Y_{t}(v_{1}, x)v_{2}=Y_{V}(v_{1}, x)v_{2}+t\Psi(v_{1}, x)v_{2}$$
for $v_{1}, v_{2}\in V$ and $(V, Y_{t}, \one)$ is a family of 
grading-restricted vertex algebras up to the first order in $t$.

The identity property for $(V, Y_{t}, \one)$ up to the first order in $t$ gives
$$Y_{V}(\one, x)v+t\Psi(\one, x)v=v+O(t^{2})$$
for $v\in V$. So we obtain
\begin{equation}\label{ident-yt}
\Psi(\one, x)v=0
\end{equation}
for $v\in V$. The creation property for $(V, Y_{t}, \one)$ 
up to the first order in $t$ gives
$$\lim_{x\to 0}(Y_{V}(v, x)+t\Psi(v, x))\one=v+O(t^{2})$$
for $v\in V$. Then we have
\begin{equation}\label{yt-one}
\lim_{x\to 0}\Psi(v, x)\one=0
\end{equation}
for $v\in V$.

The 
duality property up to the first order in $t$ can be written explicitly
as follows: For $v_{1}, v_{2}, v_{3}\in V$ and 
$v'\in V'$, 
\begin{eqnarray}
&\langle v', (Y_{V}(v_{1}, z_{1})\Psi(v_{2}, z_{2})
+\Psi(v_{1}, z_{1})Y_{V}(v_{2}, z_{2}))v_{3}\rangle&\label{prod1}\\
&\langle v', (Y_{V}(v_{2}, z_{2})\Psi(v_{1}, z_{1})
+\Psi(v_{2}, z_{2})Y_{V}(v_{1}, z_{1}))v_{3}\rangle&\label{prod2}\\
&\langle v', (Y_{V}(\Psi(v_{1}, z_{1}-z_{2})v_{2}, z_{2})
+\Psi(Y_{V}(v_{1}, z_{1}-z_{2})v_{2}, z_{2}))v_{3}\rangle&\label{iter}
\end{eqnarray}
are absolutely convergent 
in the regions $|z_{1}|>|z_{2}|>0$, $|z_{2}|>|z_{1}|>0$ and
$|z_{2}|>|z_{1}-z_{2}|>0$, respectively, to a common rational function 
in $z_{1}$ and $z_{2}$ with the only possible poles at $z_{1}, z_{2}=0$ and 
$z_{1}=z_{2}$.

Let 
\begin{eqnarray*}
Y_{V\oplus V}: (V\oplus V)\otimes (V\oplus V)&\to& (V\oplus V)[[x, x^{-1}]]\nn
(u_{1}, v_{1})\otimes (u_{2}, v_{2})&\mapsto& Y_{V\oplus V}((u_{1}, v_{1}), x)
(u_{2}, v_{2})
\end{eqnarray*}
be given by
\begin{eqnarray}\label{def-yvv}
\lefteqn{Y_{V\oplus V}((u_{1}, v_{1}), x)
(u_{2}, v_{2})}\nn
&&=(Y_{V}(u_{1}, x)u_{2}, Y_{V}(u_{1}, x)v_{2}+Y_{V}(v_{1}, x)u_{2}
+\Psi(u_{1}, x)u_{2})
\end{eqnarray}
for $u_{1}, u_{2}, v_{1}, v_{2}\in V$. 
By (\ref{def-yvv}) and (\ref{ident-yt}),
\begin{eqnarray*}
Y_{V\oplus V}((\one, 0), x)(u, v)
&=&(Y_{V}(\one, x)u, Y_{V}(\one, x)v+Y_{V}(0, x)u
+\Psi(\one, x)u)\nn
&=&(u, v)
\end{eqnarray*}
for $u, v\in V$, that is, $(V\oplus V, Y_{V\oplus V}, (\one, 0))$ 
has the 
identity property. 
By (\ref{def-yvv}) and (\ref{yt-one}),
\begin{eqnarray*}
\lefteqn{\lim_{x\to 0}
Y_{V\oplus V}((u, v), x)(\one, 0)}\nn
&&=(\lim_{x\to 0}Y_{V}(u, x)\one, \lim_{x\to 0}
Y_{V}(u, x)0+\lim_{x\to 0}Y_{V}(v, x)\one
+\lim_{x\to 0}\Psi(u, x)\one)\nn
&&=(u, v)
\end{eqnarray*}
for $u, v\in V$, that is, $(V\oplus V, Y_{V\oplus V}, (\one, 0))$ 
has the creation property. 

By (\ref{def-yvv}),
we have 
\begin{eqnarray}\label{yvv-prod1}
\lefteqn{\langle (u', v'), Y_{V\oplus V}((u_{1}, v_{1}), z_{1})
Y_{V\oplus V}((u_{2}, v_{2}), z_{2})(u_{3}, v_{3})\rangle}\nn
&&=\langle (u', v'), Y_{V\oplus V}((u_{1}, v_{1}), z_{1})\cdot\nn
&&\quad\quad\quad\quad\cdot
(Y_{V}(u_{2}, z_{2})u_{3}, Y_{V}(u_{2}, z_{2})v_{3}+Y_{V}(v_{2}, z_{2})u_{3}
+\Psi(u_{2}, z_{2})u_{3})\rangle\nn
&&=\langle u', Y_{V}(u_{1}, z_{1})Y_{V}(u_{2}, z_{2})u_{3}\rangle
+\langle v', Y_{V}(u_{1}, z_{1})
Y_{V}(u_{2}, z_{2})v_{3}\rangle\nn
&&\quad +\langle v', Y_{V}(u_{1}, z_{1})
Y_{V}(v_{2}, z_{2})u_{3}\rangle
+\langle v', Y_{V}(u_{1}, z_{1})
\Psi(u_{2}, z_{2})u_{3}\rangle\nn
&&\quad
+\langle v', Y_{V}(v_{1}, z_{1})Y_{V}(u_{2}, z_{2})u_{3}\rangle
+\langle v', \Psi(u_{1}, z_{1})Y_{V}(u_{2}, z_{2})u_{3}\rangle.
\end{eqnarray}
By the properties of $V$ and the absolute convergence of (\ref{prod1}), 
we see that the left-hand side of 
(\ref{yvv-prod1}) is absolutely convergent when $|z_{1}|>|z_{2}|>0$.
Similarly, by (\ref{def-yvv}), we have
\begin{eqnarray}\label{yvv-prod2}
\lefteqn{\langle (u', v'), Y_{V\oplus V}((u_{2}, v_{2}), z_{2})
Y_{V\oplus V}((u_{1}, v_{1}), z_{1})(u_{3}, v_{3})\rangle}\nn
&&=\langle u', Y_{V}(u_{2}, z_{2})Y_{V}(u_{1}, z_{1})u_{3}\rangle
+\langle v', Y_{V}(u_{2}, z_{2})
Y_{V}(u_{1}, z_{1})v_{3}\rangle\nn
&&\quad +\langle v', Y_{V}(u_{2}, z_{2})
Y_{V}(v_{1}, z_{1})u_{3}\rangle
+\langle v', Y_{V}(u_{2}, z_{2})
\Psi(u_{1}, z_{1})u_{3}\rangle\nn
&&\quad
+\langle v', Y_{V}(v_{2}, z_{2})Y_{V}(u_{1}, z_{1})u_{3}\rangle
+\langle v', \Psi(u_{2}, z_{2})Y_{V}(u_{1}, z_{1})u_{3}\rangle
\end{eqnarray}
and the left-hand side of 
(\ref{yvv-prod2}) is absolutely convergent when $|z_{2}|>|z_{1}|>0$.
Moreover, since (\ref{prod1}) and (\ref{prod2}) converges absolutely
when $|z_{1}|>|z_{2}|>0$ and when $|z_{2}|>|z_{1}|>0$, respectively,
to a common rational function with the only possible poles at 
$z_{1}, z_{2}, z_{1}-z_{2}=0$, 
the left-hand side of 
(\ref{yvv-prod1}) and left-hand side of 
(\ref{yvv-prod2}) also converges absolutely
when $|z_{1}|>|z_{2}|>0$ and when $|z_{2}|>|z_{1}|>0$, respectively,
to a common rational function with the only possible pole at $z_{1}-z_{2}=0$.
By (\ref{def-yvv}) again,
we have 
\begin{eqnarray}\label{yvv-iter}
\lefteqn{\langle (u', v'), 
Y_{V\oplus V}(Y_{V\oplus V}((u_{1}, v_{1}), z_{1}-z_{2})
(u_{2}, v_{2}), z_{2})(u_{3}, v_{3})\rangle}\nn
&&=\langle (u', v'), Y_{V\oplus V}(
(Y_{V}(u_{1}, z_{1}-z_{2})u_{2}, \nn
&&\quad\quad\quad\quad\quad\quad\quad\quad Y_{V}(u_{1}, z_{1}-z_{2})v_{2}
+Y_{V}(v_{1}, z_{1}-z_{2})u_{2}\nn
&&\quad\quad\quad\quad\quad\quad\quad\quad\quad\quad\quad\quad\quad\quad
+\Psi(u_{1}, z_{1}-z_{2})u_{2}), z_{2})(u_{3}, v_{3})\rangle\nn
&&=\langle u', Y_{V}(Y_{V}(u_{1}, z_{1}-z_{2})u_{2}, z_{2})u_{3}\rangle
+\langle v', Y_{V}(Y_{V}(u_{1}, z_{1}-z_{2})u_{2}, z_{2})v_{3}\rangle\nn
&&\quad +\langle v', Y_{V}(Y_{V}(u_{1}, z_{1}-z_{2})v_{2}, z_{2})u_{3}\rangle
+\langle v', Y_{V}(Y_{V}(v_{1}, z_{1}-z_{2})u_{2}, z_{2})u_{3}\rangle\nn
&&\quad +\langle v', Y_{V}(\Psi(u_{1}, z_{1}-z_{2})u_{2},
z_{2})u_{3}\rangle
+\langle v', \Psi(Y_{V}(u_{1}, z_{1}-z_{2})u_{2}, z_{2})u_{3}\rangle.\nn
\end{eqnarray}
By the properties of $V$ and the absolute convergence of (\ref{iter})
and (\ref{yvv-iter}), we see that the left-hand side of 
(\ref{yvv-iter}) is absolutely convergent when $|z_{2}|>|z_{1}-z_{2}|>0$.
Moreover, since (\ref{prod1}) and (\ref{iter}) converges absolutely
when $|z_{1}|>|z_{2}|>0$ and when $|z_{2}|>|z_{1}-z_{2}|>0$, respectively,
to a common rational function with the only possible poles at 
$z_{1}, z_{2}, z_{1}-z_{2}=0$, 
the left-hand side of 
(\ref{yvv-prod1}) and left-hand side of 
(\ref{yvv-iter}) also converges absolutely
when $|z_{1}|>|z_{2}|>0$ and when $|z_{2}|>|z_{1}-z_{2}|>0$, respectively,
to a common rational function with the only possible poles at 
$z_{1}, z_{2}, z_{1}-z_{2}=0$. So $(V\oplus V, Y_{V\oplus V}, (\one, 0))$ 
has the duality property. 

Note that the $L(-1)$-derivative property is in fact a consequence of the 
other axioms for vertex algebras. Thus $(V\oplus V, Y_{V\oplus V}, (\one, 0))$ 
is a grading-restricted vertex algebra. 

By definition, 
\begin{eqnarray*}
\lefteqn{p_{1}(Y_{V\oplus V}((u_{1}, v_{1}), x)(u_{2}, v_{2}))}\nn
&&=p_{1}(Y_{V}(u_{1}, x)u_{2}, Y_{V}(u_{1}, x)v_{2}+Y_{V}(v_{1}, x)u_{2}
+\Psi(u_{1}, x)u_{2})\nn
&&=Y_{V}(u_{1}, x)u_{2}\nn
&&=Y_{V}(p_{1}(u_{1}, v_{1}), x)p_{1}(u_{2}, v_{2})
\end{eqnarray*}
for $u_{1}, u_{2}, v_{1}, v_{2}\in V$. Also 
$$\ker p_{1}=0\oplus V$$
and 
$$Y_{V\oplus V}((0, v_{1}), x)(0, v_{2})=(0, 0)$$
for $v_{1}. v_{2}\in V$. So $p_{1}$ is a surjective homomorphism of 
grading-restricted vertex algebras and $\ker p_{1}$ is a square-zero 
ideal of $V\oplus V$.

We use $Y_{V\oplus V}^{V}$ to denote the vertex operator map
for $V\oplus V$ when $V\oplus V$ is viewed as a $V$-module. Then by definition,
\begin{eqnarray*}
i_{2}(Y_{V}(v_{1}, x)v_{2})
&=&(0, Y_{V}(v_{1}, x)v_{2})\nn
&=&Y_{V\oplus V}^{V}(v_{1}, x)(0, v_{2})\nn
&=&Y_{V\oplus V}^{V}(v_{1}, x)i_{2}(v_{2})
\end{eqnarray*}
for $v_{1}, v_{2}\in V$. So $i_{2}$ is an injective homomorphism of 
$V$-modules. Clearly, we have $i_{2}(V)=\ker p_{1}$. 
Thus $(V\oplus V, Y_{V\oplus V}, p_{1}, i_{2})$ is a 
square-zero extension of $V$ by $V$. 

Conversely, let $(V\oplus V, Y_{V\oplus V}, p_{1}, i_{2})$ be a 
square-zero extension of $V$ by $V$. Then there exists 
\begin{eqnarray*}
\Psi: V\otimes V&\to& V((x))\nn
v_{1}\otimes v_{2}&\to &\Psi(v_{1}, x)v_{2}
\end{eqnarray*}
such that 
$$Y_{V\oplus V}((u_{1}, 0), x)(u_{2}, 0)=(Y_{V}(u_{1}, x)u_{2},
\Psi(u_{1}, x)u_{2})$$
for $u_{1}, u_{2}\in V$. The identity property and 
the creation property of the 
grading-restricted vertex algebra
$(V\oplus V, Y_{V\oplus V}, (\one, 0))$ give
(\ref{ident-yt}) and (\ref{yt-one}).
The duality property for $(V\oplus V, Y_{V\oplus V}, (\one, 0))$
gives (\ref{prod1}), (\ref{prod2}) and (\ref{iter}).

For $t\in \C$, define
$$Y_{t}(v_{1}, x)v_{2}=Y_{V}(v_{1}, x)v_{2}+t\Psi(v_{1}, x)v_{2}$$
for $v_{1}, v_{2}\in V$. Then (\ref{ident-yt}) and (\ref{yt-one})
imply that $Y_{t}$ satisfies the identity property and the creation 
property up to the first order in $t$ and 
(\ref{prod1}), (\ref{prod2}) and (\ref{iter})
imply that $Y_{t}$ satisfies the duality property up to the first order in
$t$. 
Thus $(V, Y_{t}, \one)$ is a grading-restricted 
vertex algebras up to the first order in $t$, that is,
$Y_{t}$ is a first-order deformation of $(V, Y_{V}, \one)$.

Now we prove that two first-order deformations of $V$
are equivalent if and only if the corresponding 
square-zero extensions of $V$ by $V$ are equivalent.

Consider two equivalent first-order deformations of $V$ given by 
$Y_{t}^{(1)}: V\otimes V\to V((x))$ and $Y_{t}^{(2)}: V\otimes V\to V((x))$
for $t\in \C$. Then 
there exist a family
$f_{t}: V\to V$, $t\in \C$,
of linear maps of the form $f_{t}=1_{V}+tg$ where
$g: V\to V$ is a linear map preserving the grading of $V$ such that
(\ref{equiv-def}) holds 
for $v_{1}, v_{2}\in V$. 
By definition, there exist linear maps
\begin{eqnarray*}
\Psi_{1}: V\otimes V&\to& V((x))\nn
v_{1}\otimes v_{2}&\to &\Psi_{1}(v_{1}, x)v_{2}
\end{eqnarray*}
and 
\begin{eqnarray*}
\Psi_{2}: V\otimes V&\to& V((x))\nn
v_{1}\otimes v_{2}&\to &\Psi_{2}(v_{1}, x)v_{2}
\end{eqnarray*}
such that
$Y_{t}^{(1)}=Y_{V}+t\Psi_{1}$ and $Y_{t}^{(2)}=Y_{V}+t\Psi_{2}$. By 
(\ref{equiv-def}), we have 
\begin{eqnarray}\label{psi-g}
\lefteqn{\Psi_{1}(v_{1}, x)v_{2}-\Psi_{2}(v_{1}, x)v_{2}}\nn
&&=-g(Y_{V}(v_{1}, x)v_{2})+Y_{V}(g(v_{1}), x)v_{2}
+Y_{V}(v_{1}, x)g(v_{2})
\end{eqnarray}
for $v_{1}, v_{2}\in V$.

Let $(V\oplus V, Y_{V\oplus V}^{(1)}, p_{1}, i_{2})$ and 
$(V\oplus V, Y_{V\oplus V}^{(2)}, p_{1}, i_{2})$ be the 
square-zero extensions of $V$ by $V$ constructed from $Y_{t}^{(1)}$ and 
$Y_{t}^{(2)}$. Let $h: V\oplus V\to V \oplus V$ 
be defined by 
$$h(v_{1}, v_{2})=(v_{1}, v_{2}+g(v_{1}))$$
for $v_{1}, v_{2}\in V$. Clearly, $h$ is a linear isomorphism. 
For $u_{1}, u_{2}, v_{1}, v_{2}\in V$, by definition and (\ref{psi-g}),
\begin{eqnarray*}
\lefteqn{h(Y^{(1)}_{V\oplus V}((u_{1}, v_{1}), x)(u_{2}, v_{2}))}\nn
&&=h(Y_{V}(u_{1}, x)u_{2}, Y_{V}(u_{1}, x)v_{2}+Y_{V}(v_{1}, x)u_{2}
+\Psi_{1}(u_{1}, x)u_{2})\nn
&&=(Y_{V}(u_{1}, x)u_{2}, Y_{V}(u_{1}, x)v_{2}+Y_{V}(v_{1}, x)u_{2}\nn
&&\quad\quad \quad\quad \quad\quad \quad\quad \quad
+\Psi_{1}(u_{1}, x)u_{2}+g(Y_{V}(u_{1}, x)u_{2}))\nn
&&=(Y_{V}(u_{1}, x)u_{2}, Y_{V}(u_{1}, x)v_{2}+Y_{V}(v_{1}, x)u_{2}\nn
&&\quad\quad \quad\quad \quad\quad \quad\quad \quad
+\Psi_{2}(u_{1}, x)u_{2}+Y_{V}(g(u_{1}), x)u_{2}
+Y_{V}(u_{1}, x)g(u_{2})))\nn
&&=(Y_{V}(u_{1}, x)u_{2}, Y_{V}(u_{1}, x)(v_{2}+g(u_{2}))\nn
&&\quad\quad \quad\quad \quad\quad \quad\quad \quad
+Y_{V}((v_{1}+g(u_{1})), x)u_{2}
+\Psi_{2}(u_{1}, x)u_{2})\nn
&&=Y^{(2)}_{V\oplus V}(h(u_{1}, v_{1}), x)h(u_{2}, v_{2}).
\end{eqnarray*}
So $h$ is in fact an isomorphism from the algebra
$(V\oplus V, Y_{V\oplus V}^{(1)}, (\one, 0))$ to the algebra
$(V\oplus V, Y_{V\oplus V}^{(2)}, (\one, 0))$.
Now it is clear that the following diagram is 
commutative:
$$\begin{CD}
0@>>> V @>>i_{2}> V\oplus V @>>p_{1}>V @>>>0\\
@. @V1_{W}VV @VhVV @VV1_{V}V\\
0@>>> V @>>i_{2}> V\oplus V @>>p_{1}>V @>>>0,
\end{CD}$$
So these two first order deformations are equivalent.

Conversely, let $(V\oplus V, Y_{V\oplus V}^{(1)}, p_{1}, i_{2})$ and 
$(V\oplus V, Y_{V\oplus V}^{(2)}, p_{1}, i_{2})$ be two 
equivalent 
square-zero extensions of $V$ by $V$. 
Let $\Psi_{1}, \Psi_{2}: V\otimes V
\to V((x))$ be given by
\begin{eqnarray*}
Y_{V\oplus V}^{(1)}((u_{1}, 0), x)(u_{2}, 0))
&=&(Y_{V}(u_{1}, x)u_{2}, \Psi_{1}(u_{1}, x)u_{2}),\\
Y_{V\oplus V}^{(2)}((u_{1}, 0), x)(u_{2}, 0))
&=&(Y_{V}(u_{1}, x)u_{2}, \Psi_{2}(u_{1}, x)u_{2})
\end{eqnarray*}
for $u_{1}, u_{2}\in V$. Then $Y_{t}^{(1)}, 
Y_{t}^{(2)}: V\otimes V
\to V((x))$ given by
\begin{eqnarray*}
Y_{t}^{(1)}(v_{1}, x)v_{2}=Y_{V}(v_{1}, x)v_{2}+t\Psi_{1}(v_{1}, x)v_{2},\\
Y_{t}^{(2)}(v_{1}, x)v_{2}=Y_{V}(v_{1}, x)v_{2}+t\Psi_{2}(v_{1}, x)v_{2}
\end{eqnarray*}
for $v_{1}, v_{2}\in V$ are first-order deformations of $(V, Y_{V}, \one)$
by the proof above. 

Let $h: V\oplus V\to V\oplus V$ be an equivalence from 
$(V\oplus V, Y_{V\oplus V}^{(1)}, p_{1}, i_{2})$ to
$(V\oplus V, Y_{V\oplus V}^{(2)}, p_{1}, i_{2})$. Then by
Lemma \ref{h-form}, there exists
a linear map $g: V\to V$ such that
$$h(v_{1}, v_{2})=(v_{1}, v_{2}+g(v_{1}))$$
for $v_{1}, v_{2}\in V$. Using the fact that $h$ is an isomorphism of 
grading-restricted vertex algebras from 
$(V\oplus V, Y_{V\oplus V}^{(1)}, (\one, 0))$ to 
$(V\oplus V, Y_{V\oplus V}^{(2)}, (\one, 0))$,
we obtain (\ref{psi-g}) which implies (\ref{equiv-def}). 
Thus the two first-order deformations $Y_{t}^{(1)}$ and $Y_{t}^{(2)}$
are equivalent. 
\epfv

\noindent {\small \sc Beijing International Center for Mathematical
Research, Peking University, Beijing 100871, China},

\vspace{1em}

\noindent {\small \sc Kavali Institute For Theoretical Physics 
China, CAS, Beijing 100190, China}

\vspace{1em}

\noindent  {\it and}

\vspace{1em}
\noindent {\small \sc Department of Mathematics, Rutgers University,
110 Frelinghuysen Rd., Piscataway, NJ 08854-8019  (permanent address)}

\vspace{1em}

\noindent {\em E-mail address}: yzhuang@math.rutgers.edu

\end{document}